%Paper: funct-an/9405003
%From: Karl H. Hofmann <hofmann@mathematik.th-darmstadt.de>
%Date: Mon, 16 May 94 7:42:31 MESZ

%%%%%%% THIS IS EPI.TEX= a plain TeX file on a memo on epis of C$^*$-algebras

%%%%%%% FORMAT MACROS %%%%%%

\magnification=\magstep1
\tolerance=300
\pretolerance=200
\hfuzz=1pt
\vfuzz=1pt

\hsize=5in
\vsize=8in

%% We print on A4 sheets:

\hoffset=210 true mm             %A4 page width
\voffset=297 true mm             %A4 page height
\advance\hoffset by -\hsize      %compute complete margin space
\divide\hoffset  by 2            %left margin
\advance\hoffset by -1 true in   %standard offset is 1 in
\advance\voffset by -\vsize      %compute complete margin space
\divide\voffset  by 2            %top margin
\advance\voffset by -1 true in   %standard offset is 1 in

\parindent=35pt
\mathsurround=1pt
\parskip=1pt plus .25pt minus .25pt
\normallineskiplimit=.99pt

%%%%%%% FONT MACROS %%%%%%%

\font\bfone=cmbx10 scaled\magstep1
\def\qed{{\unskip\nobreak\hfil\penalty50\hskip .001pt \hbox{}\nobreak\hfil
          \vrule height 1.2ex width 1.1ex depth -.1ex
           \parfillskip=0pt\finalhyphendemerits=0\medbreak}\rm}

%%%%%%% TEXT FOLLOWS %%%%%%%%

\vglue1truecm
\centerline{\bfone Epimorphisms of C$^*$-algebras are surjective}
\bigskip
\centerline{\bf K. H. Hofmann and  K.-H. Neeb}
\bigskip

In this note we observe that epimorphisms in the category of $C^*$-algebras
are surjective [2]. To show this claim it clearly suffices
to prove that epimorphically embedded subalgebras cannot be proper, and this is
what we shall verify.\medskip\noindent
\bf Definition 1.\quad\rm We say that a subalgebra $B$ of the
$C^*$-algebra $A$ is \it epimorphically embedded \rm
if the inclusion morphism $B\to A$ is an epimorphism in the category of
C$^*$-algebras.
\qed\medskip\noindent
\bf Lemma 2.\quad\it Let $B$ be an epimorphically embedded subalgebra of $A$
and $(\pi, {\cal H})$ a representation of $A$ on a Hilbert space ${\cal H}$.
Suppose that $v\in {\cal H}$ is a cyclic vector for $A$. Then $v$ is also a
cyclic vector for  $B$.\smallskip\noindent
\bf Proof.\quad\rm
Set ${\cal K}= \overline{\pi(B)v}$. We must show ${\cal H}\subseteq{\cal K}$.
Now ${\cal K}$ is a closed subspace
of ${\cal H}$ which is invariant under $B$. Let $P$ denote the orthogonal
projection onto ${\cal K}$ and set $U={\bf1}-2P$. Then $U$ is unitary
and $U(v_1+v_2)=-v_1+v_2$ for $v_1\in{\cal K}$ and
$v_2\in{\cal K}^\bot$.

We define a new representation $\pi'\colon A\to L({\cal H})$ by
$\pi'(x)=U\pi(x)U^{-1}$. Since the subalgebra $B$ commutes with $P$ and
thus with $U$, it follows that $\pi|B=\pi'|B$. The
assumption that $B$ is epimorphically embedded now yields $\pi=\pi'$.
Hence $U$ commutes with $\pi(A)$.
Thus ${\cal K}$, the eigenspace of $U$ for the eigenvalue $-1$,
is invariant under $A$. Therefore $v\in {\cal K}$
implies ${\cal H}=\overline{\pi(A)v}\subseteq{\cal K}$.\qed\medskip\noindent
\bf Lemma 3.\quad\it
Suppose that $B$ is an epimorphically embedded subalgebra of $A$ and
$\phi$ and $\phi'$
are two positive linear forms on $A$ with $\phi|B=\phi'|B$.
 Then $\phi=\phi'$.
\smallskip\noindent
\bf Proof.\quad\rm
Let $(\pi,{\cal H})$ and $(\pi',{\cal H}')$ denote
the corresponding representations of $A$ with the cyclic vectors
$v\in {\cal H}$ and $v'\in {\cal H}'$, respectively, satisfying
$$ \phi(a)=\langle\pi(a)v,v\rangle\quad\hbox{and}\quad
\phi'(a)=\langle\pi'(a)v',v'\rangle$$
for all $a\in A$ (cf. [1], 2.4.4, p.~32). We first
claim that for $w=\pi(b)v$, $b\in B$ the
element $\pi'(b)v'$ depends on $w$ but not on the choice of $b\in B$.
Indeed suppose that $w=\pi(b')v$. Then set $c=b'-b\in B$. Now $\pi(c)v=0$ and
$\phi|B=\phi'|B$ implies
$$0=\langle\pi(c)v,\pi(c)v\rangle=\phi(c^*c)=\phi'(c^*c)
=\langle\pi'(c)v',\pi'(c)v'\rangle.$$
Now that the first claim is proved
we define unambiguously a linear map
$$\Psi_0\colon\pi(B)v\to\pi'(B)v'\quad\hbox{by}\quad
  \Psi_0(\pi(b)v)=\pi'(b)v'.$$
Using $\phi|B=\phi'|B$ again we note
$$\eqalign{\langle\Psi_0\big(\pi(b)v\big),\Psi_0\big(\pi(b')v\big)\rangle
&=\langle\pi'(b)v',\pi'(b')v'\rangle\cr
&=\phi'(b'^*b)=\phi(b'b^*)=\langle\pi(b)v,\pi(b')v\rangle.\cr}$$
Hence  $\Psi_0$ is a linear  isometry.
But by Lemma 2, the subspaces $\pi(B)v$ of ${\cal H}$ and $\pi'(B)v'$
of ${\cal H}'$, respectively, are dense.
Thus $\Psi_0$ extends uniquely to a linear
isometry $\Psi\colon{\cal H}\to{\cal H}'$.

This provides us with two representations of $A$ on ${\cal H}'$, namely,
$\pi'$ and $\rho$ given by $\rho(a)=\Psi\circ\pi(a)\circ\Psi^{-1}$.
This is equivalent to $\rho(a)\circ\Psi=\Psi\circ\pi(a)$.
We claim that $\pi'|B=\rho$. To see this, let $b,\,b'\in B$. Then
$$\eqalign{\rho(b)\big(\pi'(b')v'\big)
&=\rho(b)\Psi\big(\pi(b')v\big)
 =\Psi\big(\pi(b)\pi(b')v\big)
 =\Psi\big(\pi(bb')v\big)
 =\pi'(bb')v'\cr
&=\pi'(b)\big(\pi'(b')v'\big).\cr}$$
Now the density of $\pi'(B)v$ in ${\cal H}'$ (established in Lemma 2)
implies our claim.

Since $B$ is epimorphically embedded in $A$, it follows that
$\pi' = \rho$. As a consequence,
$$\eqalign{\phi'(a)
&=\langle\pi'(a)v',v'\rangle=\langle\rho(a)v',v'\rangle\cr
 =\langle\Psi\pi(a)\Psi^{-1}v',v'\rangle\cr
&=\langle\pi(a)\Psi^{-1}v',\Psi^{-1}v'\rangle=\langle\pi(a)v,v\rangle\cr
&=\phi(a).\cr}$$
This shows that $\phi=\phi'$ as asserted.\qed\medskip
For a C$^*$ algebra $A$ let $A_h$ denote the closed vector subspace
of hermitian elements of $A$.

\noindent
\bf Lemma 4.\quad\it
Let $A$ be a C$^*$-algebra and $B$ a proper closed
subalgebra of $A$. Then there exist two different positive linear forms
$\phi$ and $\phi'$ of $A$ such that $\phi|B = \phi'|B$.\smallskip\noindent
\bf Proof.\quad\rm
(Compare [1], Lemma 11.3.2 and its proof (p.~228).)
Since $B\ne A$ we have an element $a\in A_h\setminus B_h$. By the Hahn-Banach
Theorem we find a continuous linear functional $f_h$ on $A_h$ with
$f_h(B_h)=\{0\}$ and $f_h(a)=1$. Now the definition $f(x+iy)=f_h(x)+if_h(y)$
(for $x,\,y\in A_h$) extends $f_h$  uniquely to a continuous
hermitean linear form $f$ on $A$ with $f(B)=\{0\}$ and $f(a)=1$. By
[1] 2.6.4 (p.~40) there are two positive forms $\phi$ and $\phi'$ with
$f=\phi-\phi'$. Then $\phi|B=\phi'|B$ and $\phi(a)\ne\phi'(a)$.\qed\medskip
\noindent
\bf Theorem 5.\quad \it
 If a {\rm C}$^*$-subalgebra $B$ of a {\rm C}$^*$-algebra $A$ is
epimorphically embedded then $B=A$.\smallskip\noindent
\bf Proof.\quad \rm The assertion is a direct consequence of Lemmas 3 and 4.
\qed \medskip\noindent
[1]\quad  Dixmier, J., Les C$^*$-alg\`ebres et leurs repr\'esentations,
Gauthier-Villars, Paris, 1964$^1$, 1969$^2$.
\smallskip\noindent
[2]\quad Pestov, V.G., Are epimorphisms of C$^*$-algebras onto?,
Internet newsgroup sci.math.research posting, 6 April 1994, article
No. 2179.
\smallskip\noindent
[3] ---, New features of noncommutativity in the presence of topology, in:
Abstracts of Talks, 1994 New Zealand Mathematical Colloquium,
University of Waikato, Hamilton, May 1994.
\bye